\documentclass{amsart}
\usepackage[dvips]{graphics}

\usepackage{graphicx}
\input{epsf.tex}
\newtheorem*{thm}{Theorem}
\newtheorem*{thm1}{Theorem 1}
\newtheorem*{thm2}{Theorem 2}

\newtheorem*{corollary}{Corollary}
\newtheorem*{lem2}{Lemma 2}
\newtheorem*{lem1}{Lemma 1}
\newtheorem*{lem3}{Lemma 3}
\newtheorem*{lem4}{Lemma 4}

\newtheorem*{definition}{Definition}
\newtheorem*{remark}{Remark}
\newtheorem*{proposition}{Proposition}

\newcommand {\bd}{\partial}
\newcommand{\leqslant}{\leq}
\newcommand{\geqslant}{\geq}
\newcommand{\G}{\gamma}

\newcommand{\ra}{\rightarrow}
\newcommand{\la}{\leftarrow}

\numberwithin{equation}{section}

\begin{document}

\title{Asymptotic dimension of finitely presented groups}

\thanks{Research co-funded by the European Social Fund and National
Resources-EPEAEK II-PYTHAGORAS}

\author{Thanos Gentimis}
\address{Mathematics Department, University of Athens,
Athens 157 84, Greece } \email{gentimisth@math.uoa.gr}


\date{\today}


\begin{abstract}
We prove that if a finitely presented group is one-ended then its
asymptotic dimension is bigger than $1$. It follows that a finitely
presented group of asymptotic dimension $1$ is virtually free.
\end{abstract}

\maketitle

\section{Introduction}

The notion of asymptotic dimension of a metric space was introduced
by Gromov in \cite {GRO2}. It is a large scale analog of topological
dimension and it is invariant by quasi-isometries. This notion has
proved relevant in the context of Novikov's higher signature
conjecture and it was investigated further by other people (see
\cite {Yu}, \cite {Dr},  \cite {RO}).

In this paper we show the following:
\begin{thm1} If $G$ is a one-ended finitely presented group then $G$ has asymptotic dimension
greater or equal to 2.
\end{thm1}

Also we deduce as corollary:
\begin{thm2} If $G$ is a finitely presented group with $asdim\,G =
1$ then $G$ is virtually free.
\end{thm2}

For finitely generated groups the statement above doesn't hold. We
give a counter-example at the end.

After we completed this work T.Januszkiewicz brought to our
attention his joint paper with J.Swiatkowski (\cite {J-S}) where
they proved the same results independently.

Finally I would like to thank pr. P.Papasoglu for his help and
guidance during the writing of this paper.

\section{Preliminaries}

\textbf{Metric Spaces}. Let $(X,d)$ be a metric space. If $A,B$ are 
subsets of $X$ we set $d(A,B)=\inf \{d(a,b):a\in A, b\in B \}$. A 
\textit{path} in $X$ is a map $\G : I \ra X $ where $I$ is an 
interval in $\mathbb R$. A path $\G $ joins two points $x$ and $y$ 
in $X$ if $I = [a,b]$, $\G(a) = x$ and $\G(b) = y$. The path $\G$ is 
called an \textit{infinite ray} starting from $x_0$ if 
$I=[0,\infty)$ and $\G(0)=x_0$. A \textit{geodesic}, a 
\textit{geodesic ray}, or a \textit{geodesic segment} in X, is an 
isometry $\G : I \ra X$ where $I$ is $\mathbb R$ or $[0,\infty)$ or 
a closed interval in $\mathbb R$. We use the terms geodesic, 
geodesic ray etc. for the images of $\G$ without discrimination. On 
a path-connected space $X$ given two points $x,y$ we define the 
\textit{path metric} to be $\rho(x,y) = inf \{length(p)\}$ where the 
infimum is taken over all paths $p$ that join $x$ and $y$.  A space 
is called a \textit{geodesic metric space} if for every $x$, $y$ in 
$X$ there exists a geodesic segment which joins them. In a geodesic 
space the path metric is indeed a metric. A geodesic metric space 
$X$ is said to be \textit{one-ended} if for every bounded $K$, $X-K$ 
has exactly one unbounded connected component. We say that $X$ is 
\textit{uniformly one-ended}, if for every $n\in \mathbb R^+$ there 
is an $m\in \mathbb R^+$ such that for every $K \subset X$ with 
$diam(K)<n$, $X-K$ has exactly one connected component of diameter 
bigger than $ m$. 

\

\textbf{Groups}. The \textit{Cayley graph} of $G$ with respect to a 
generating set $S$ is the 1-dimensional complex having a vertex for 
each element of $G$ and an edge joining vertex $x$ to vertex $xs$ 
for every vertex $x$ and every $s\in S$  . The Cayley graph has a 
natural metric which makes it a geodesic metric space, where each 
edge has length 1 (see  \cite{RO}). In fact any connected graph can 
be made geodesic metric space in the same way.

We will use the same letter $G$ for both the group and its Cayley
graph as a metric space. We also use \textit{van - Kampen
Diagrams} (see \cite{L-S} chapter V pg 236 -240 and \cite{Ger}). A
van Kampen diagram $\mathfrak D$ for a
 word $w$ in $S$ representing the identity element of $G$, is a finite, planar, contractible,
combinatorial 2-complex; its 1-cells are directed and labeled by
generators and the boundary labels of each of its 2-cells are cyclic
conjugates of relators or inverse relators. Further the boundary
label for $\mathfrak D$ is $w$ when read (by convention
anticlockwise) from a base point in $\bd \mathfrak D$. We recall
here that a word $w$ represents the identity element of $G$ if and
only if the path in the Cayley graph labeled by $w$ is closed.

 There is a
natural map $f$ from the 1-skeleton $\mathfrak D^{(1)}$ of
$\mathfrak D$ to the Cayley graph of $G$. $f$ sends the base point
to a vertex $v$ of the Cayley graph and edges of $\mathfrak
D^{(1)}$ to edges of the Cayley graph with the same label.
Obviously $f$ is determined by the image of the base point, $v$.
 $f$ is
not necessarily injective. If we consider $\mathfrak D^{(1)}$ as a
geodesic metric space giving each edge length 1 then
$d(x,y)\geqslant d(f(x),f(y))$ for every $x,y$ in $\mathfrak
D^{(1)}$.

We say that a group $G$ is \textit{virtually free} if there exists a 
finite index subgroup $H$ of $G$ which is a free group.

A group $G$ is called free-by-finite if it contains a normal
subgroup $N$ of finite index such that $N$ is free.

Generally if $X$ is a property of groups then we say that $G$ is
virtually $X$ if $\exists\, H<_f G$ and $H$ has property $X$. We say
that $G$ is $X£$-by-finite if $\exists\,N\lhd_f G$ normal such that
$N$ has property $X$.

Obviously if $G$ is $X$-by-finite then $G$ is virtually $X$. The
converse also holds if the property $X$ is inherited to subgroups.

Now since if $G$ is free every subgroup of $G$ is free we have that
$G$ is virtually free if and only if $G$ is free by-finite.

\

\textbf{Asymptotic Dimension}.
 A metric space $Y$ is said to be $d$ - disconnected or that
it has dimension 0 on the $d$ - scale if  there exist $B_i \subset
Y$ such that: $$ Y \,\, = \bigcup_{i\in I}
 B_i $$ with $\sup \{ diam B_i,\,\, i \in I \} \leqslant D <
 \infty $, d($B_i,B_j$) $\geqslant d$ $\forall \,\, i \neq j.$
\newline
\begin{definition}\textit{(Asymptotic Dimension 1)} We say that a space $X$ has
asymptotic dimension $n$ if $n$ is the minimal number such that for
every $d > 0$
 we have: \newline $X =\bigcup X_k$ for $k$ = $0,1,2, ... ,n$ and all $X_k$ are $d$-disconnected. We then write
  \newline $asdim\,X =
 n$.
\end{definition}
 We say that a covering  $\{B_i\}$ of $X$ has \textit{d - multiplicity} $\leq k$, if and only if every 
$d$-ball $B(x,d)$ in $X$ meets no more than $k$ sets $B_i$ of the 
covering. The covering has \textit{d - multiplicity} $k$ if 
\textit{d - multiplicity} $\leq k$ is true and \textit{d - 
multiplicity} $\leq k-1$ is false. A covering $\{ B_i \}, \,\,i \in 
I$ is \textit{D - bounded}, if $diam(B_i) \leqslant D $ for all 
$i\in I.$
\begin{definition} \textit{(Asymptotic Dimension 2)} We say that a
space $X$ has $asdim X = n$, if $n$ is the minimal number such
that $\forall \,d > 0$ there exists a $D$-bounded covering of $X$
with $d$ - multiplicity $\leqslant n+1$.
\end{definition}

 The two definitions used here are the first two definitions Gromov gave in his paper
 \cite{GRO2}.  It is not difficult to see that the two definitions are
 equivalent.

\section{Main Theorem}

Before we get to the main theorem we will prove two lemmas that we
will need below. \begin{lem1} Let $G$ be a finitely generated,
infinite group then asdim $G>0$. \end{lem1}

\begin{proof} Let asdim $G = 0$ and fix  $d>0$. Then according to
the first definition we have that $G = X_1$ were $X_1 = \bigcup B_i$
with:
\begin{enumerate}
\item $diam B_i \leqslant D$, $\forall \,i \in I$ \item
$d(B_i,B_j) \geqslant d$, $\forall \,i\neq j.$
\end{enumerate}
That means that $G$ is $d$-disconnected. Since G is a connected
graph it follows that we can not have two distinct $B_i$'s. So
$G\subset B_1$ which means that $G$  is $D$-bounded. But since $G$
is finitely generated we have immediately that $G$ is finite which
is a contradiction.
\end{proof}
\begin{lem2} If $G$ is an one-ended finitely presented group then
$G$ contains a bi-infinite geodesic.
\end{lem2}
\begin{proof} Take any $n\,\in\,\mathbb N$. Then define \[C_n=\{\text{geodesic paths starting from e of length
n}\}\] We note that:
\begin{itemize}
\item[a)] $C_n \neq \emptyset$ for all $n\,\in\,\mathbb N$.
\item[b)] $C_n$ is a finite set for all $n\,\in\,\mathbb N$.
\end{itemize}

Consider a map \[\pi_{n_1}^n: C_n\ra C_{n-1}\] which takes every geodesic path
of length $n$ from $G_n$ and cuts off the last edge of that path.
Then obviously what is left is also a geodesic and now the length is
$n-1$. Thus it is contained in $C_{n-1}$. So for every
$n\,\in\,\mathbb N$, $\pi_{n-1}^n$ is a natural well defined map. 
$\newline\newline$ So now consider the inverse limit sequence
\[ \{e\}\la C_1\la C_2 \la ... \la C_n \la C_{n+1} \la ... \]
with bonding maps $\pi^n_{n-1}$. By using the well known statement 
that the inverse limit of compact spaces is compact and hence 
nonempty, we obtain that 
\[F=\lim_{\la}\{C_i,\pi^{i}_{i-1}\}\ne\emptyset\] 
Now let $r\,\in\,F$. Obviously $r$ is an infinite geodesic ray 
starting from $e$. Now lets fix a geodesic segment $\xi$ in $C_n$ 
and denote by $Y_\xi$ all the infinite geodesics in F starting with 
$\xi$. Then obviously $Y_{\xi}$ is the inverse image 
$(\pi^\infty_n)^{-1}(\xi)$ where $\pi^{\infty}_n:F\to C_n$ is the 
natural projection. Thus if $Y_\xi$ is non empty it is compact.  
$\newline\newline$ Now define:
\[S_n = \{g\,\in\,G\,\, \text{such that} \,\, d(g,e)=n\},\] 
\[L_n=\{\text{all the geodesic paths starting from}\,\, S_n\,\, 
\text{and ending at}\,\, e\},\] \[E_n = \{r: r \,\,\text{are 
geodesic rays},\,\, r(0)\,\in\,K_n, \,\,e\,\in\,r_i\}\] Obviously 
$E_n$ is not empty since let $r'$ one of the geodesic rays we 
defined starting from $e$. Let $g=r'(n)$ be the vertex of $r'$ in 
$S_n$. Define $g^{-1}\cdot r'$. Since multiplying with an element 
doesn't change the respective distances $g^{-1}\cdot r'$ is also a 
geodesic starting from $g^{-1}\cdot e=g^{-1}$ which is obviously in 
$S_n$ and passing through $g^{-1}\cdot g=e$. Thus, $E_n\neq 
\emptyset$. 

Now fix an $\eta$ in $L_n$ and its corresponding $h\in S_n$. Define 
the inverse limit \[\{h\}\la h\cdot C_1\la h\cdot C_2\la ... \la 
h\cdot C_n \la ...\] The above inverse limit gives us all the 
geodesic rays starting from $h$ i.e. $h\cdot F$. 

The set $Y'_\eta$ of all geodesic rays starting with the path $\eta$ 
being homeomorphic to $Y_{h^{-1}\cdot\eta}$, is a compact  by the 
above remark (if it is not empty of course!).

Then by the definition we get that:
\[E_n=\bigcup_{\eta\,\in\,L_n}Y'_\eta\]
Since $E_n$ is non-empty, it is compact as a finite union of compact 
sets. Thus the inverse limit:\[E=\lim_{\la}E_n\] with the 
restrictions as the bonding maps is compact and nonempty. Now let 
$k$ in $E$. Then $k$ is a path that passes through $e$ and $k$ is 
our bi-infinite geodesic line. This concludes the lemma.
\end{proof}

\begin{thm1} If $G$ is an one-ended finitely presented group then
asdim $G \geqslant 2$.
\end{thm1}

\begin{proof}
By lemma $1$ and since $G$ is one-ended we have that $G$ is infinite
and thus asdim $G>0$ We will show that asdim $G \neq 1$. Let's
suppose that asdim $G = 1$.

 Let $M$ = $\max \{|r_i|,i=1,2,...,n$: $r_i$ relation of $G\}$, where
 $|r| = length$ of the word $r$. We fix $d>100M+100$. Since asdim $G = 1$  there is a covering $\mathbb B = \{B_i\}$
 with:
$$G=\bigcup_{i \in I} B_i$$ and $diam B_i <D$, $\forall
\,i\,\in \,I$, such that every ball $B(x,d)$ intersects at most $2$
sets of the covering $\mathbb B$. We may assume without loss of
generality that if $r$ is a path in the Cayley graph labeled by a
relator $r_i$ then $r$ is contained in some $B_j\in \mathbb B$.

 Since $G$ is one-ended we have that $G$ has a
bi-infinite geodesic $S$(Lemma $2$). Let $N =
max\{100D^{100},300M\}$. Choose an $x_0 \in S $ and consider the
ball $B(x_0,N)$ which separates the geodesic $S$ into two geodesic
rays $S_1$ and $S_2$. Since $G$ is one-ended there is an $x$ in
$S_1$, a $y$ in $S_2$ and a path $p$ with $p(0) = x$ and $p(t) = y$
such that $p\bigcap B(x_0,N)= \emptyset$.

 We denote by $[x,y]$ the part of the geodesic $S$, that connects $x$ and $y$.
Obviously $length([x,y])\geqslant 2N$. We denote by $w$ the path
that corresponds to $[x,y]\bigcup p$ \newline We have then that
$$length(w) = length([x,y]) + length(p) \Rightarrow length(w)> 200D$$
 So in order to cover the path $w$ we
need at least $3$ sets of the covering $\{B_i\}$. \newline We
consider now the van-Kampen diagram $\mathfrak D$ that corresponds
to the path $w$ and the function $f$ from $\mathfrak D^{(1)}$ to
the Cayley graph $G$. So $f(\partial \mathfrak D)=w$. For
notational convenience we label vertices and edges of $\partial
\mathfrak D$ in the same way as $w$. So for example we denote  the
vertex on $\partial \mathfrak D$ which is mapped to $x_0\in w$ by
$f$ also by $x_0$.

 Let $B$ be a set of the covering that
intersects $[x,y]$. We consider $f^{-1}(B)$. Let $C(B)$ be the
union of all 2-cells of $\mathfrak D$ which have the property that
their boundary is contained in $f^{-1}(B)$. Let $U$ be the
collection of all such sets $B$ with the following property:
 For some connected component, $K$, of $C(B)$,  $[x,y]\cap K$ is
 contained in an interval $[a,b]$ with $a,b\in K$ such that
 $x_0\in [a,b]$.
 Let $d(K)=d(a,b)$ for such a component and let $d(B)$ be the
 maximal value of all $d(K)$ for $K$ component of $C(B)$ such that
 $x_0\in [a,b]$.
 Let $B_1$ be
a set in $U$ for which $d(B_1)$ is maximal. Let $K_1$ be the
connected component of $C(B_1)$ for which $d(K_1)=d(B_1)$. Let's
say that $K_1\cap [x,y]$ is contained in $[a_1,b_1]$ with
$a_1,b_1\in K_1$. Let $e$ be the edge of $[x,y]$ adjacent to $a_1$
which does not lie in $[a,b]$. If $r$ is the 2-cell containing $e$
there is some $B_2\in \mathbb B$ such that $C(B_2)$ contains $r$.

%
%
%

Let $C$ be the subset of $\mathfrak D$ which contains all 2-cells
with boundary contained in $f^{-1} (B_1\bigcup B_2)$. Since
$d(x_0,p)\geq N$, $C$ does not intersect $p$. Thus $\mathfrak
D-C\neq \emptyset$. Let $K$ be the connected component of $C$
which contains $x_0$.

Let $P = \overline{\mathfrak D-K}\bigcap \overline{K}$. $P$
 is connected since $K$ is connected. Each edge of $P$ is
 contained in two 2-cells. One of these 2-cells lies in $f^{-1} (B_1\bigcup B_2)$
 and one does not lie in this set. It is not possible that all
 edges of $P$ are contained in a 2-cell of $f^{-1} (B_2)$. Indeed in
 this case we would have $d(B_2)>d(B_1)$, which is impossible.
 Since $r$ is contained in $C(B_2)$ some edge of $P$ is not
 contained in a 2-cell of $f^{-1} (B_1)$. It follows that there
 are 2 adjacent edges $e_1,e_2$ in $P$ such that one of them is
 contained in a 2-cell of $f^{-1} (B_1)$ and the other in a 2-cell
 of $f^{-1} (B_2)$. If $c$ is the 2-cell that contains $e_1$ and is
 not in $f^{-1} (B_1\bigcup B_2)$ then $c$ lies in a set
 $f^{-1}(B_3)$ with $B_3\in \mathbb B$ and $B_3\ne B_1,B_2$. The
 edges $e_1,e_2$ and the 2-cell $c$ have a vertex $v$ in common.
 So $v\in f^{-1} (B_1\cap B_2\cap B_3)$. It follows that $B_1\cap
 B_2\cap B_3 \ne \emptyset $, a contradiction.

%
This concludes the proof.
\end{proof}

\begin{remark} The result above holds also for uniformly one-ended
simply connected simplicial complexes. So if $X$ is a uniformly
one-ended simply connected simplicial complex then $asdim\,X\geq
2$.
\end{remark}


 Of course the result does not hold for one-ended simply
connected simplicial complexes, a half-line gives a
counterexample. We remark finally that if a Cayley graph is one
ended then it is uniformly one ended.

We note that the following theorem holds:
\begin{thm}(Dunwoody-Stallings \cite{D})
If $G$ is a finitely presented group then G is the fundamental group
of a  graph of groups such that all the edge groups are finite and
all the vertex groups are $0$ or $1$ ended.
\end{thm}
Also it is known that:
\begin{lem3} If all the vertex groups are $0$-ended (i.e. finite) then G is
virtually free (see \cite {S}, page 120, prop.11).
\end{lem3}

Furthermore it is not difficult to prove the following lemma (see
\cite{GRO1}):
\begin{lem4} If $H<G$ and $H$ is finitely generated then $asdim \,G
\geqslant asdim \,H$.
\end{lem4}

Using the lemma above and theorem $1$ we have the stronger result:

\begin{thm2} If $G$ is a finitely presented group with asdim $G =
1$ then $G$ is virtually free.
\end{thm2}
\begin{proof} Let $G$ be a finitely presented group with asdim $G = 1$. Let $\Gamma$ be the
graph of groups of the Dunwoody-Stallings theorem. If a vertex group
$H$ is one-ended then from the theorem $1$, asdim $H \geqslant 2$.
But $H<G$ which means that asdim $G\geqslant 2$ which is a
contradiction. So all vertex groups are $0$-ended. It follows that
$G$ is virtually free.
\end{proof}

Now we give an example of a finitely generated group which is not
finitely presented, not virtually free and that has asymptotic
dimension $1$. Namely:

\begin{proposition} Let $G=\mathbb Z_2\wr_{r}\mathbb Z$ be the restricted wreath
product of $\mathbb Z_2$ and $\mathbb Z$. Then asdim $G = 1$ and $G$
is not virtually free.
\end{proposition}

\begin{proof} Since $G=\mathbb Z_2\wr_{r}\mathbb Z$ there exists a short exact sequence:
\[0\rightarrow(\oplus_{\mathbb Z} \mathbb Z_2)\rightarrow G \rightarrow \mathbb
Z\rightarrow 0\] By the Hurewicz type formula (see \cite{DS}) we
have: \[asdim G \leq asdim\mathbb Z + asdim (\oplus_{\mathbb Z}
\mathbb Z_2)\]

Since every finitely generated subgroup $F$ of $(\oplus_{\mathbb Z}
\mathbb Z_2)$ is finite, we have that all $F$ of that type have
asymptotic dimension $0$. Following the definition of asymptotic
dimension for arbitrary discrete groups found in \cite{DS} we get:
\[asdim (\oplus_{\mathbb Z} \mathbb Z_2) = \sup \{\text{asdim} F|
F<G, \text{finitely generated} \}=0\] Another way to get the same
result is by using the following corollary found in \cite{J}
\begin{corollary} Let G be a countable abelian group. Then asdim $G=0$ if and only if $G$ is
torsion.\end{corollary} Obviously $\oplus_{\mathbb Z} \mathbb Z_2$
is abelian and torsion so asdim $\oplus_{\mathbb Z} \mathbb Z_2=0$.
Thus we get:
\[ \text{asdim } G \leq 1+0=1\] Since $G$ is a finitely generated infinite
group (see \cite{X} for the description of the generators), by lemma
$1$ we have that asdim $G>0$. So asdim $G = 1$.

We will prove that $G$ is not virtually free. Let $G$ be virtually
free. Then $G$ is free-by-finite. So there exists a normal subgroup
$N$ of $G$ such that $N$ is free and the index $|G:N|$ is finite.
Recall the exact sequence:

\[0\rightarrow(\oplus_{\mathbb Z} \mathbb Z_2)\rightarrow G
\rightarrow \mathbb Z\rightarrow 0\] and lets denote the image of
$\oplus_{\mathbb Z} \mathbb Z_2$ in $G$ to be $H$ through the one to
one mapping $f$. Then obviously $H<G$. Define $q:G\ra \frac{G}{N}$
and restrict $q$ to $H$. Since $H$ is infinite and $\frac{G}{N}$ is
finite then $H\bigcap kerq\neq \emptyset$. But $kerq=N$ thus
$H\bigcap N \neq \emptyset$ which is a contradiction because if
$a\,\in\,H$ $a$ is torsion. Thus $a\,\notin N$ since $N$ is free.
\end{proof}
\bibliographystyle{amsplain}

\end{document}